\documentclass[12pt]{article}   	
\usepackage[margin=2.5cm]{geometry}                		
\usepackage{graphicx}
\usepackage{amsmath,amssymb,mathrsfs,amsthm}
\usepackage{soul}
\usepackage{mathtools}
\usepackage{subcaption}
\theoremstyle{definition}
\usepackage{float}
\usepackage[flushleft]{threeparttable}
\theoremstyle{theorem}

\usepackage{color}

\newcommand{\lrp}[1]{\left(#1\right)}
\newcommand{\lrb}[1]{\left[#1\right]}
\newcommand{\lrc}[1]{\left\{#1\right\}}
\renewcommand{\exp}[1]{\textrm{exp}\lrc{#1}}
\renewcommand{\log}[1]{\textrm{log}\lrp{#1}}

\renewcommand{\d}{\text{d}}

\usepackage{abstract,lipsum}

\newtheorem{thm}{Theorem}

\title{Heavy Tailed Horseshoe Priors}
\author{Andrew Womack, Zikun Yang}
\date{}							
\begin{document}
\maketitle
\begin{abstract}
Locally adaptive shrinkage in the Bayesian framework is achieved through the use of local-global prior distributions
that model both the global level of sparsity as well as individual shrinkage parameters for mean
structure parameters. The most popular of these models is the Horseshoe prior and its variants due to their
spike and slab behavior involving an asymptote at the origin and heavy tails. In this article, we present an
alternative Horseshoe prior that exhibits both a sharper asymptote at the origin as well as heavier tails, which
we call the Heavy-tailed Horseshoe prior. We prove that mixing on the shape parameters provides improved
spike and slab behavior as well as better reconstruction properties than other Horseshoe variants. A simulation study is provided to show the advantage of the heavy-tailed Horseshoe in terms of absolute error to both the truth mean structure as well as the oracle.

\end{abstract}
\section{Introduction}
We are interested in the typical normal mean estimation problem under and an assumption of sparsity of signal. Suppose that the data ${\pmb y}$ is a $n\times 1$ vector, where each data point $y_i$ is normally distributed with mean $\mu+\phi_i$ and variance $\sigma^2$.  We propose the following model:
\begin{equation}
\label{eq::prior}
\begin{array}{rclrcl}
y_i|\mu,\sigma^2,\phi_i&\stackrel{ind}{\sim}& \mathcal{N}(\mu+\phi_i,\sigma^2);&\phi_i|\gamma_i,\sigma^2,Z &\stackrel{ind}{\sim}& \mathcal{N}(0,\frac{\sigma^2}{\gamma_i Z});\\
\mu,\sigma^2,Z &\sim& \pi(\mu,\sigma^2,Z);&\gamma_i &\stackrel{iid}{\sim}& \mathcal{G}(\gamma_i),
\end{array}
\end{equation}
where $\gamma_i$ is the locally individual-level adaptive parameter and $Z$ is the global adaptive parameter. The assumption of sparsity is that most of the $\phi_i$ are exactly $0$.\\

The prior distribution of the global parameters ($\sigma^2$, $\mu$, and $Z$) is not our primary interest in this paper. For $\sigma^2$ and $\mu$, a typical Normal-Inverse-Gamma family with weak prior information has been proved to be flexible and is the default prior for these parameters. As for $Z$, its prior should provide a good multiplicity control and work as an indication of the overall sparsity. The choice of the prior distribution of $Z$ is based mostly on the data set at hand and the corresponding sampling scheme. Common default choices are a gamma distribution or a beta distribution of the second kind. Our main focus in this paper is the prior distribution $G$ on the $\gamma_i$ and its effects on sparsity.

\subsection{The Horseshoe and Horeseshoe+ Priors}
At face value, the ability to locally adapt comes directly from  the individual-level random effect parameter $\phi_i$. Ideally, one would expect that  $\phi_i$ is estimated close to zero if the data point is close to the global mean $\mu$ or provide a good estimate of the distance between $\mu$ and $y_i$ if the data point is not close to $\mu$. Define the local shrinkage profile parameter $\tau_i=\frac{\gamma_i}{1+\gamma_i}\in\lrp{0,1}$ and assume $\sigma^2=1,~\mu=0$, and $Z=1$, then the posterior mean of $\phi_i$ can be represented by the following  equation
\begin{equation}
\label{eq::post_phi}
\text{\pmb {E}}\lrb{\phi_i|y}=y_i-\text{\pmb E}\lrb{\tau_i|y_i}y_i.
\end{equation}
This is showing that the amount of posterior shrinkage is controlled by the  $\tau_i$, where $\tau_i\to 1$ meaning total shrinkage and $\tau_i\to 0$ meaning no shrinkage. Hence, the focus is on the $\gamma_i$ or $\tau_i$ for controlling the $\phi_i$ for this type of model. \\

In \cite{carvalho2010horseshoe}, the authors give a thorough review of the previous priors of $\tau_i$ and propose the famous Horseshoe shaped Beta$\lrp{\frac{1}{2},\frac{1}{2}}$ for $\tau_i$. Immediately, the Horseshoe prior drew a significant amount of attention due to its bounded influence on the estimation of $\phi_i$ and its noise control property. A further reason for the popularity that the Horseshoe prior receives is the relatively easy sampling scheme through the following hierarchy: 
\begin{equation}
\label{eq::HS_hier}
\begin{array}{rcl}
\gamma_i|\omega_i&\sim&\text{Gamma}\lrp{0.5,\omega_i}\\
\omega_i&\sim&\text{Gamma}\lrp{0.5,1}.
\end{array}
\end{equation}
It is easy to observe that this hierarchy results in a complete conjugate sampling procedure for the most of the parameters of the Horseshoe model, hence implementing the Horseshoe prior   is very straight forward through Gibbs sampling.\\

The Horseshoe prior is more of a good inspiration for the global-local adaptive mechanism than it is a practical estimation and prediction tool. Oftentimes, it fails to provide enough shrinkage to perform well in the ultra-sparse situation. In \cite{bhadra2017horseshoe+}, the authors propose adding two more layers of the Gamma distribution into the hierarchy in \eqref{eq::HS_hier} in order to attain a better result than the original Horseshoe prior. They refer this prior as the Horseshoe plus (HS+). The marginal densities of the HS(+) are 
\begin{equation}
\label{table::prior_compare}
\begin{array}{rcl}
\pi_{HS}\lrp{\gamma_i}&=&\frac{\gamma_i^{-0.5}}{\pi(\gamma_i+1)}\\
\pi_{HS+}\lrp{\gamma_i}&=&\frac{\gamma_i^{-0.5}}{\pi(\gamma_i-1)}\frac{\log{\gamma_i}}{\pi}.
\end{array}
\end{equation}
Notice that  at the tail area both marginals are dominated by $\gamma_i^{-1.5}$ term and that near the origin both marginals are dominated by $\gamma_i^{-0.5}$ term. The marginal prior from the HS+ is only different from the marginal prior from the HS up to a slowly varying function of $\gamma_i$. This suggests a minimal improvement for the HS+ over the HS. There is clearly room for improvement in creating priors for the $\gamma_i$ that are dramatically different than those for the HS.

%
%
\subsection{The Heavy Tail Horseshoe Prior}

We propose to include another individual-level parameter $p_i$, which could be referred as the local decision parameter. This parameter serves as the shape parameter in the Gamma hierarchy in \label{eq::HS_hier}, which is modified to
\begin{equation}
\label{eq::hths_hier}
\begin{array}{rcl}
\gamma_i|p_i,\omega_i&\sim&\text{Gamma}\lrp{p_i,\omega_i}\\
\omega_i|p_i&\sim&\text{Gamma}\lrp{1-p_i,1}\\
p_i&\sim&\pi(p_i).
\end{array}
\end{equation}
A first, and simplest, model assumes $p_i\sim \text{Uniform}\lrp{0,1}$. This provides a closed form of  the marginal density of $\gamma_i$, which will be discussed in this section. We further investigate the other possible choices of the priors of  $p_i$ later in the paper. The relationship between $p_i$, $\gamma_i$, and $\omega_i$ can be revealed from \eqref{eq::hths_hier}. Note that $p_i\to1$ causes $\omega_i$ to be smaller and $\gamma_i$ to be larger. At same time, the smaller $\omega_i$ reinforces a larger $\gamma_i$. Similar logic holds for $p_i\to 0$. This decision-reinforcing mechanism characterizes the advantage of our model compared to the original HS prior of the HS+ prior. The effects of this decision reinforcement will be shown through the comparison of the prior marginal densities of $\gamma_i$ and $\tau_i$ as well as posterior effects.\\

 In the hierarchy \eqref{eq::hths_hier}, it is possible to mathematically integrate out $\omega_i$ and $p_i$,  providing a closed form for the marginal prior of $\gamma_i$ and $\tau_i$. The marginal density of $\gamma_i$ is
\begin{equation}
\label{eq::g_m_prior}
\pi(\gamma_i)=\frac{1}{\gamma_i\lrp{\lrp{\log{\gamma_i}}^2+\pi^2}};~\gamma_i\in\lrp{0,\infty}.
\end{equation}
It turns out that the marginal density of $\gamma_i$ is a log-Cauchy density function, and the $\frac{1}{{\lrp{\log{\gamma_i}}^2+\pi^2}}$ part barely makes the density function integratabtle on the positive real line. Hence, it is an extreme heavy tailed prior density function. We call it the heavy tailed horseshoe prior, refer as `HTHS' in this article. Similar to the HS+ hierarchical structure, an even heavier tail prior density can be achieved by introducing another two layers of Gamma distribution upon the HTHS prior hierarchy, then integrating the hierarchy out to acquire the marginal. Strikingly, we are surprised to observe that the resulting density form is 
\begin{equation}
\label{eq::g_m_prior_+}
\pi(\gamma_i)=\frac{2}{\gamma_i\lrp{\lrp{\log{\gamma_i}}^2+4\pi^2}};~\gamma_i\in\lrp{0,\infty},
\end{equation}
where it is also a log-Cauchy distribution with only a different scaling for $\log(\gamma_i)$. It will be referred as `HTHS+' in this article. Comparing to the HS(+) prior densities, the marginal densities of the HTHS(+) are both proportion to $\gamma_i^{-1}$ at the origin and at the tail up to a slowly varying function $h(\gamma_i)=\lrp{\log{\gamma_i}}^2+\pi^2$. It indicates that the HTHS pushes more probability mass toward the origin and the tail areas to have a more decisive judgement on the nature of $\phi_i$, coinciding the initial idea of adding the local decision parameter $p_i$. \\

\begin{figure}[H]
\centering
\begin{minipage}{.5\textwidth}
  \centering
  \includegraphics[width=3in]{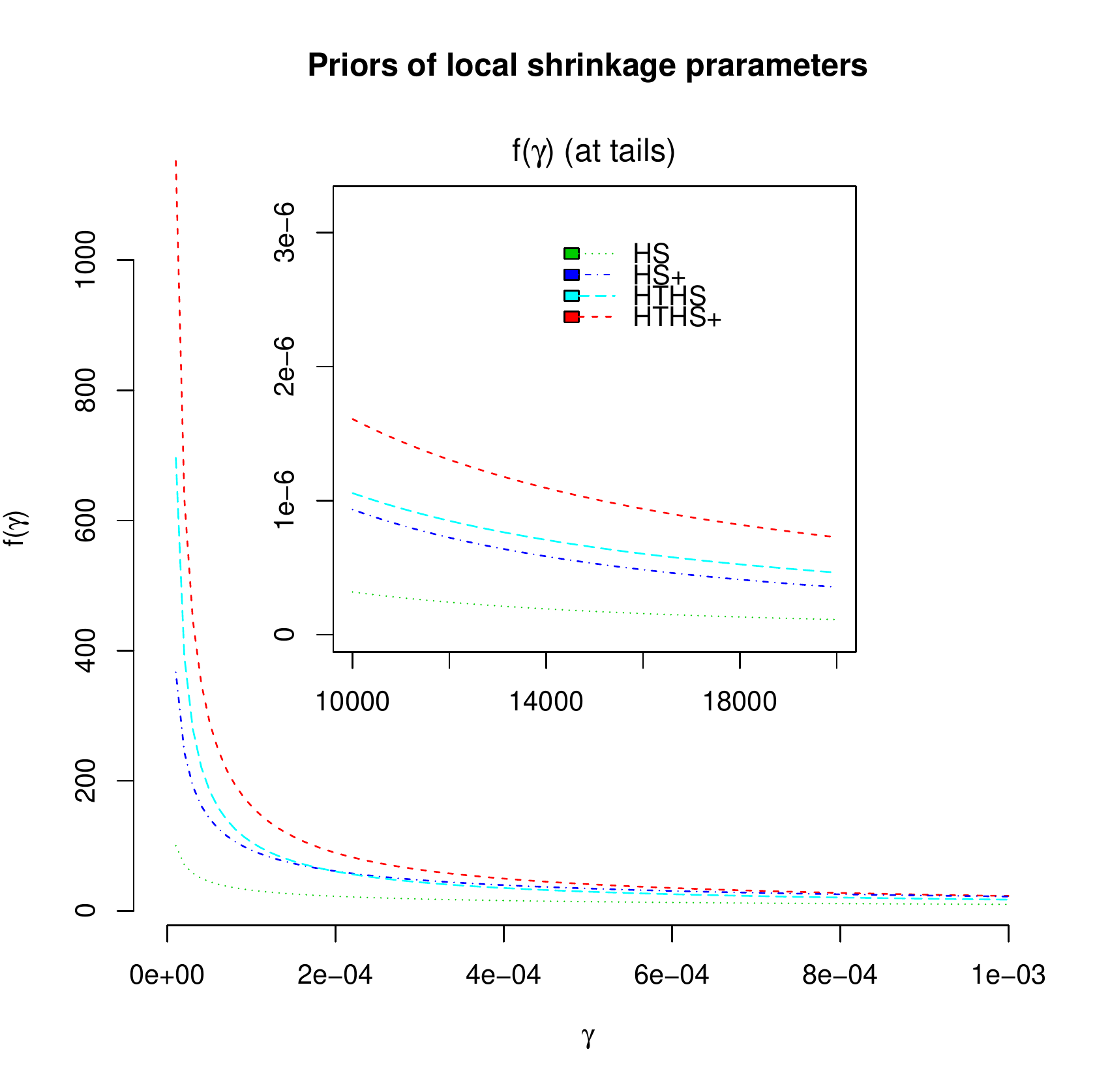}
  \caption{Prior densities of $\gamma_i$.}
   \label{fig:gamma_compare}
\end{minipage}%
\begin{minipage}{.5\textwidth}
  \centering
  \includegraphics[width=3in]{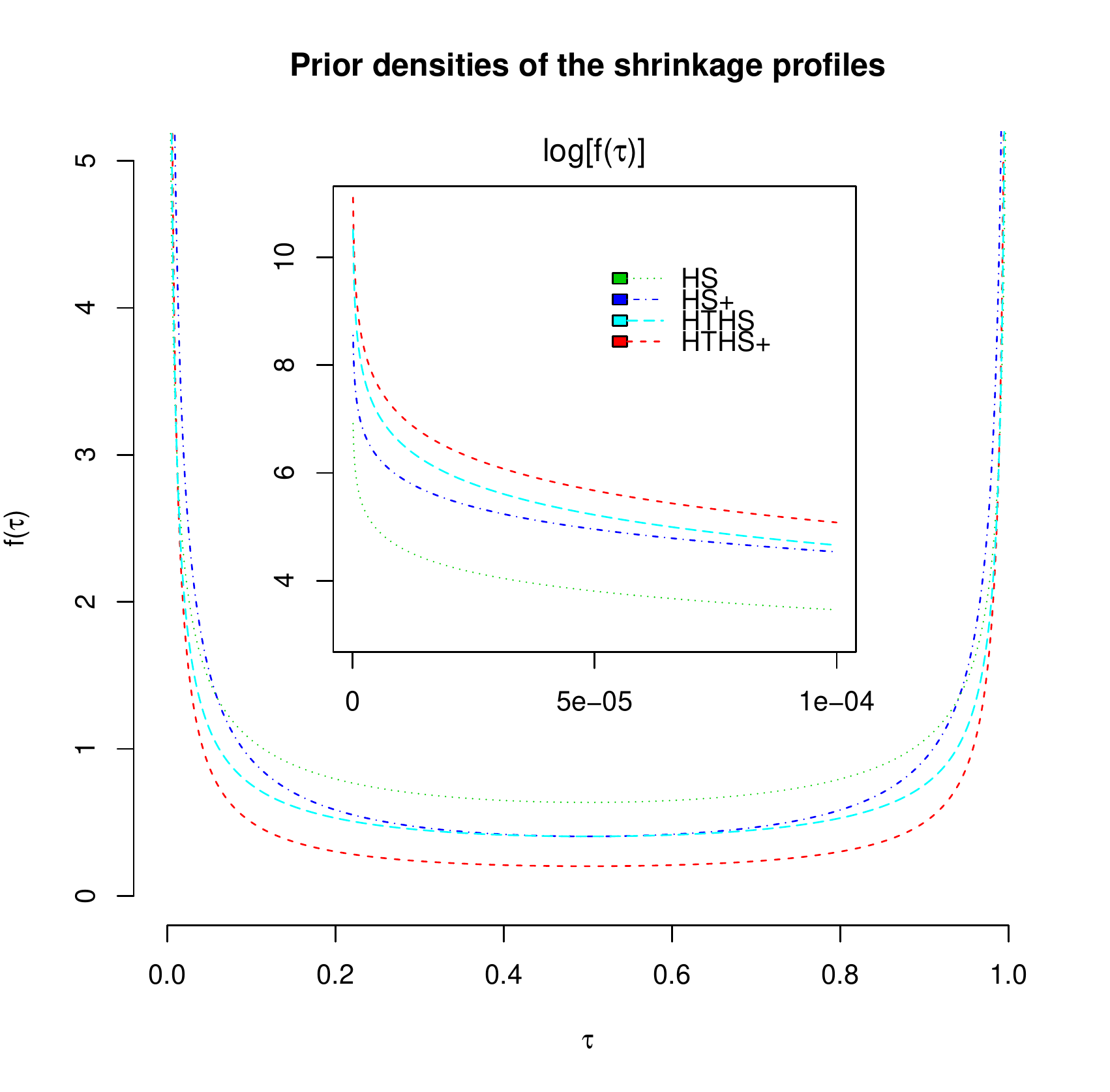}
   \caption{Prior densities of $\tau_i$. }
   \label{fig::tau_compare}
\end{minipage}
\end{figure}

We could also observe the improvement from the marginals of the shrinkage profile parameter $\tau_i$. 
\begin{equation}
\label{eq::kappa_compare}
\begin{array}{rclcrcl}
\pi_{HS}\lrp{\tau_i}&=&\frac{\lrb{\tau_i\lrp{1-\tau_i}}^{-0.5}}{\pi}&\quad&\pi_{HS}\lrp{\tau_i}&=&\frac{\lrb{\tau_i\lrp{1-\tau_i}}^{-0.5}}{\pi}
\frac{\log{\tau_i}-\log{1-\tau_i}}{\pi(2\tau_i-1)}\\
\pi_{HTHS}\lrp{\tau_i}&=& \frac{\lrb{\tau_i\lrp{1-\tau_i}}^{-1}}{\lrb{\log{\tau_i}-\log{1-\tau_i}}^{2}+\pi^2}&\quad&\pi_{HTHS+}\lrp{\tau_i}&=& \frac{2\lrb{\tau_i\lrp{1-\tau_i}}^{-1}}{\lrb{\log{\tau_i}-\log{1-\tau_i}}^{2}+4\pi^2}.
\end{array}
\end{equation}
Similar to the arguments made for $\gamma_i$, the ratios between the prior densities of the HS(+) and the HTHS(+) are $\tau_i^{-0.5}$ at the origin and $\lrp{1-\tau_i}^{-0.5}$ at the unit (up to a slowly varying function), pushing more probability mass away from the center area and avoiding ambiguous decisions. Figure \ref{fig:gamma_compare} and  Figure \ref{fig::tau_compare} shows a graphical comparison for these priors.
\section{Theoretical Properties}
This section provides statements of the theoretical advantage of adopting the HTHS over the HS prior, especially under the sparse condition.  

\subsection{Marginal Distribution for $y$}
\label{sec::margin}
The following theorem characterizes the tail behavior of the marginal likelihood of the HTHS(+).
\begin{thm}
\label{thm::likeli_margin}
Suppose $y\sim N(\phi,1)$ and that $\phi\sim HTHS(+)$ with $Z=1$.  Let $m(y)=\int N(y|\phi,1)f(\phi)\text{d}\phi$ denote the marginal density for $y$. 
Then as $y\to \infty$, the marginal density satisfies
\begin{equation}
\label{eq::tail_y}
m(y)
\stackrel{|y|\rightarrow \infty}{\asymp}
\frac{L(y^2)}{|y|},
\end{equation}
where $L$ is a slowly varying function (for every $t>0$, $L(tx)/L(x)\rightarrow 1$ as $x\rightarrow\infty$). As a consequence, we have
\begin{equation}
\label{eq::logscore_y}
\frac{\text{d}}{\text{d}y}\log{m(y})\stackrel{|y|\rightarrow \infty}{\asymp} \frac{1}{|y|},
\end{equation}
up to the (more quickly vanishing) score of the slowly varying function $L$. 
\end{thm}
It is striking to observe that both likelihood marginals of the HS(+) are proportion to $\frac{1}{y^2}$ up to slowly varying functions.  Comparing to the HS(+),  the marginal distribution of the HTHS(+) deploys more probability mass on large values as the marginal distribution is asymptotically $\frac{1}{|y|}$. This theorem completely captures the advantage of the HTHS(+) for preserving large signals. In \cite{carvalho2010horseshoe}, the authors point out that this kind of robustness for large features comes from the heavy tail of the prior density of $\phi_i$. Hence, it is not surprising to see that the tails of $\pi_{HTHS(+)}(\phi_i)$ are heavier than the tails of $\pi_{HS(+)}$, which will be proved in the next subsection. 
\subsection{Marginal of $\phi_i$}
The impact of $\gamma_i$ on the marginal density of the individual random effect parameter $\phi_i$ is one of the key points of all the locally adaptive models. However, neither of these prior densities of $\phi_i$ have an analytically closed forms due to the complexities of the marginals of $\gamma_i$. It is true that all four prior densities share basic features, such as they all have an asymptote at the origin and are in the domain of polynomial tails. Assume $\sigma^2=1$, $Z=1$, and $\mu=0$, the following theorem gives the asymptotic upper and lower bounds of the marginal density of $\phi_i$ from the HTHS(+),
\begin{thm}
\label{thm::LUbounds}
Let $a\in(0,0.5)$. The univariate marginal densities of for $\phi$ from the HS+ and HTHS+ satisfy the inequalities
\begin{equation}
\label{eq::ineqs}
\begin{array}{rcl}
C_a f(\phi|HS+)&\leq&
\frac{a2^{a-1}}{\sqrt{\pi}|\phi|^{1+2a}}\Gamma_L\left(0.5+a,\frac{\phi^2}{2}\right) + 
\frac{a2^{-a-1}}{\sqrt{\pi}|\phi|^{1-2a}}\Gamma_U\left(0.5-a,\frac{\phi^2}{2}\right)
\\
&\leq &
D_a f(\phi|HTHS+)
\end{array}
\end{equation}
where $\Gamma_L$ and $\Gamma_U$ are the lower and upper incomplete gamma functions, respectively, and $C_a$ and $D_a$ are positive numbers.
\end{thm}
By the nature of the incomplete gamma functions, Theorem \ref{thm::LUbounds} indeed proves the features of spikiness and slabbiness of the HTHS(+). \\

The lack of the tractable analytical forms of the HS(+) still prevents us from direct insights on the comparison of the tail behaviors of $\phi_i$ between the HS(+) and HTHS(+), which is crucial for the estimation robustness. It turns out that a simple integration trick could shed the light on this matter if we keep the intermediate hierarchical parameter $\omega_i$ from the whole integration procedure. We have
\begin{equation}
\label{eq::tail}
\begin{array}{rclcl}
\pi_{HS}\lrp{\phi}&\propto&\int_0^\infty \frac{1}{\omega+\frac{\phi^2}{2}}\exp{-\omega}  \d \omega&\stackrel{|\phi|\to\infty}{\asymp}&\frac{1}{\phi^2}\\
\pi_{HS+}\lrp{\phi}&\propto&\int_0^\infty  \frac{\log{\omega}-\log{\frac{\phi^2}{2}}}{\omega-\frac{\phi^2}{2}}\exp{-\omega} \d \omega&\stackrel{|\phi|\to\infty}{\asymp}&\frac{\log{|\phi|}}{\phi^2}\\
\pi_{HTHS}\lrp{\phi}&\propto&\int_0^\infty  \frac{\lrp{\frac{\phi^2}{2}+\omega+1}}{\lrb{
\lrp{\log{\frac{\phi^2}{2}+\omega}}^2+\pi^2
}}\frac{1}{
\lrp{
\frac{\phi^2}{2}+\omega
}^{1.5}
}\exp{-\omega} \d \omega&\stackrel{|\phi|\to\infty}{\asymp}&\frac{1}{|\phi|\times\log{|\phi|}^2}.
\end{array}
\end{equation}
Equation \eqref{eq::tail} directly presents the fact that the marginal of $\phi$ from the HTHS has a heavier tail than the HS(+), further  concluding the result of the marginal likelihood robustness in section \ref{sec::margin}. 
\begin{figure}[H] 
   \centering
   \includegraphics[width=5in]{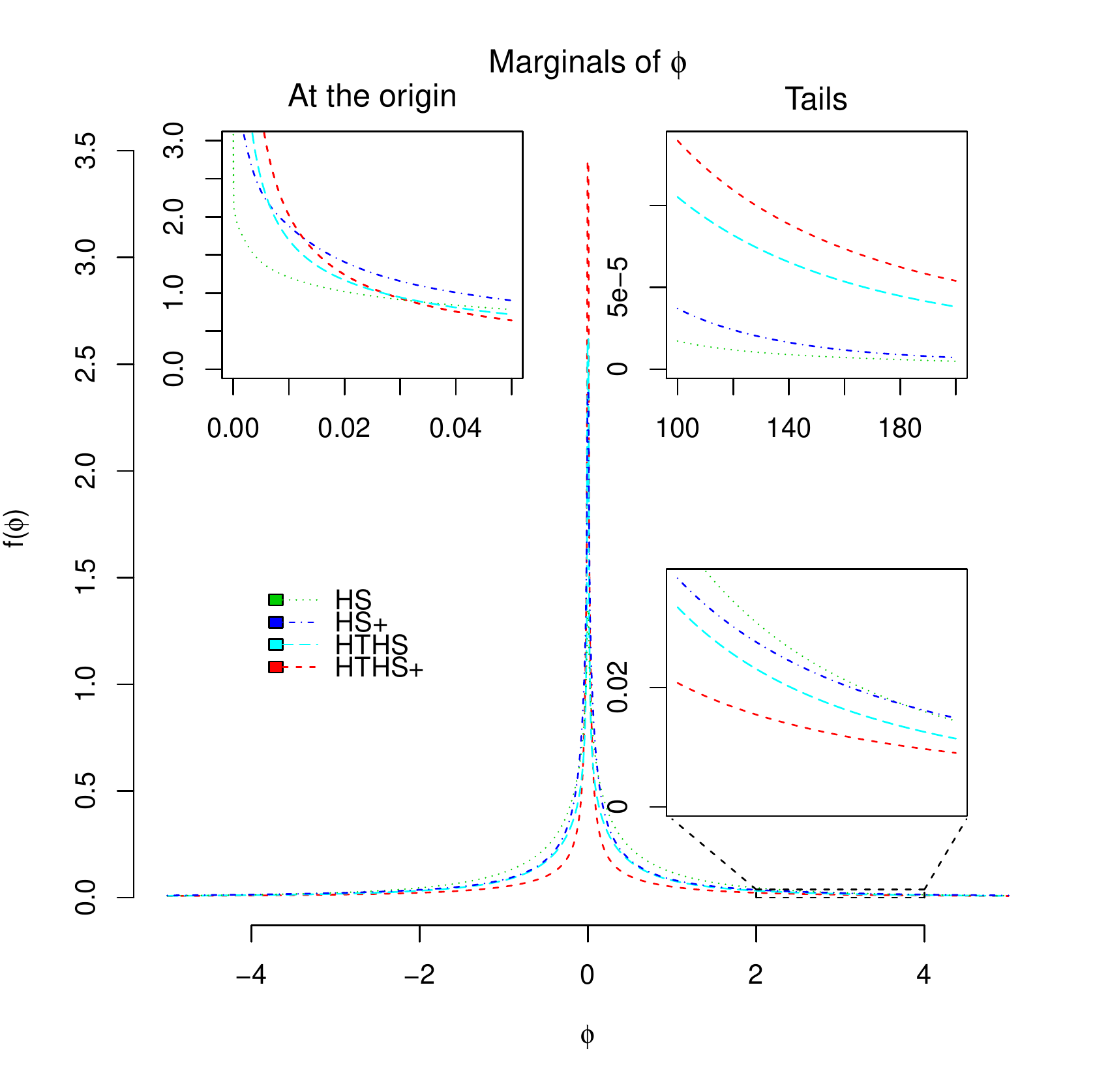} 
   \caption{The prior densities of $\phi$}
   \label{fig:phi}
\end{figure}
Figure \ref{fig:phi} shows the prior marginal densities of $\phi$ from the HS(+) and the HTHS(+) based on numerical approximation. Notice that the HTHS(+) is more spikiness at the origin than the HS(+), and this spikiness is critical to pursue a better risk bound when the true signal is indeed sparse by providing more shrinkage on the noise. The theorem in the next section gives a consideration on the risk bound. 
\subsection{K-L risk bounds}
When the true value of each $\phi_i$ is zero, all of the locally-adaptive shrinkage models are super-efficient (convergence to the true mean at a rate faster than that of the MLE). This is due to the prior asymptote at the origin. One way to measure their relative efficiencies is through the Kullback-Leibler risk bound, which characterizes the divergence between  the true sampling density and the Bayesian predictive density. Let $L(f_1,f_2)$ denote the K-L divergence of $f_2$ from $f_1$, the Ces\`{a}ro average Bayes predictive risk upper bound for true value $\phi_0$ is
\begin{equation}
\label{eq::cesaro}
R_n(\phi_0)\le \epsilon-n^{-1}\log{\pi\left(A_\epsilon\right)},
\end{equation}
where $\pi$ is the prior for $\phi$ and $A_\epsilon=\lrc{\phi:L(f(y|\phi_0),f(y|\phi))<\epsilon}$. The following theorem gives the Kullback-Leibler risk bounds of the HTHS(+) prior at the origin. 

\begin{figure}
\begin{minipage}{\linewidth}
   \centering
   \includegraphics[width=.4\linewidth]{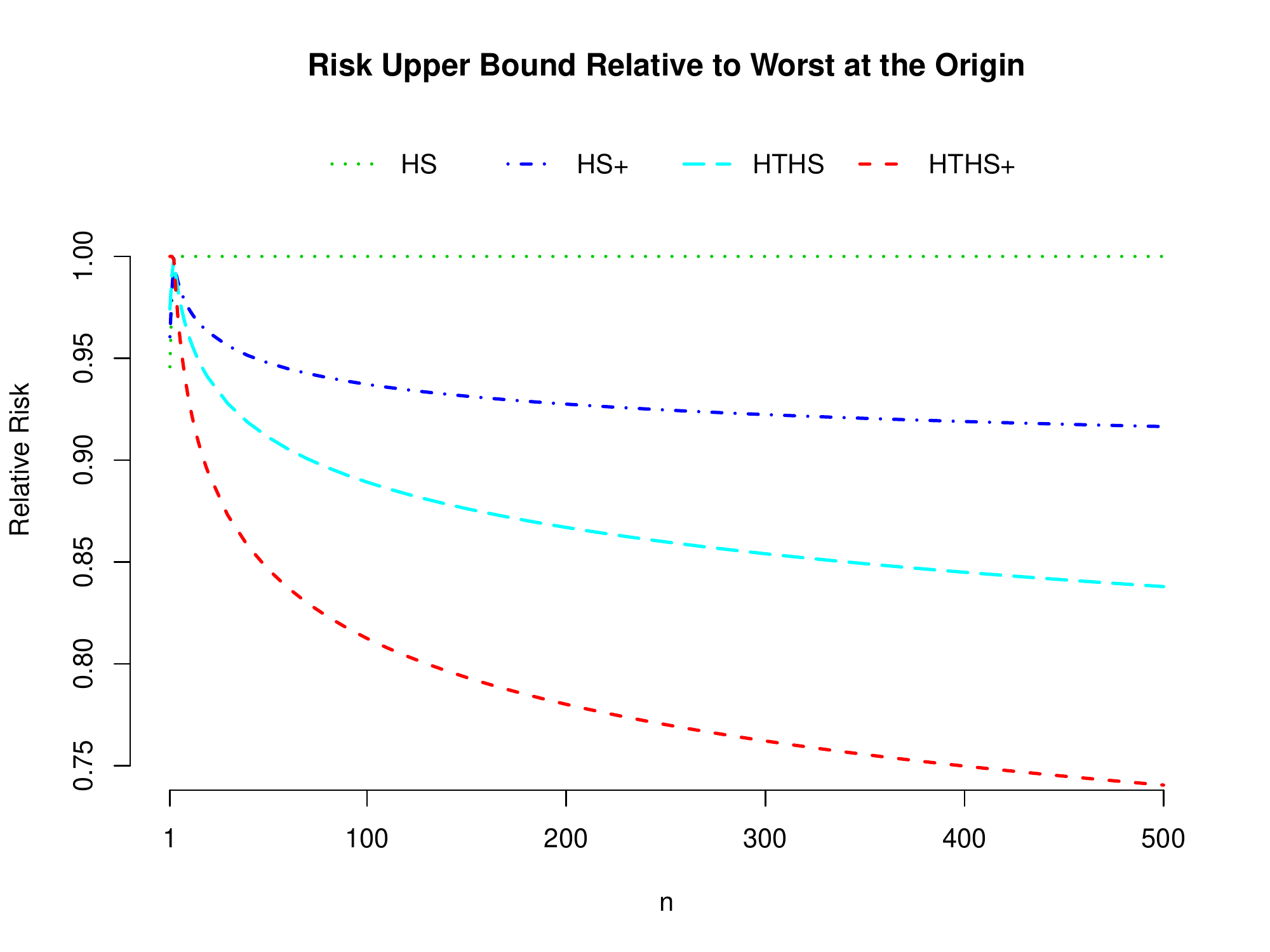} 
   \includegraphics[width=.4\linewidth]{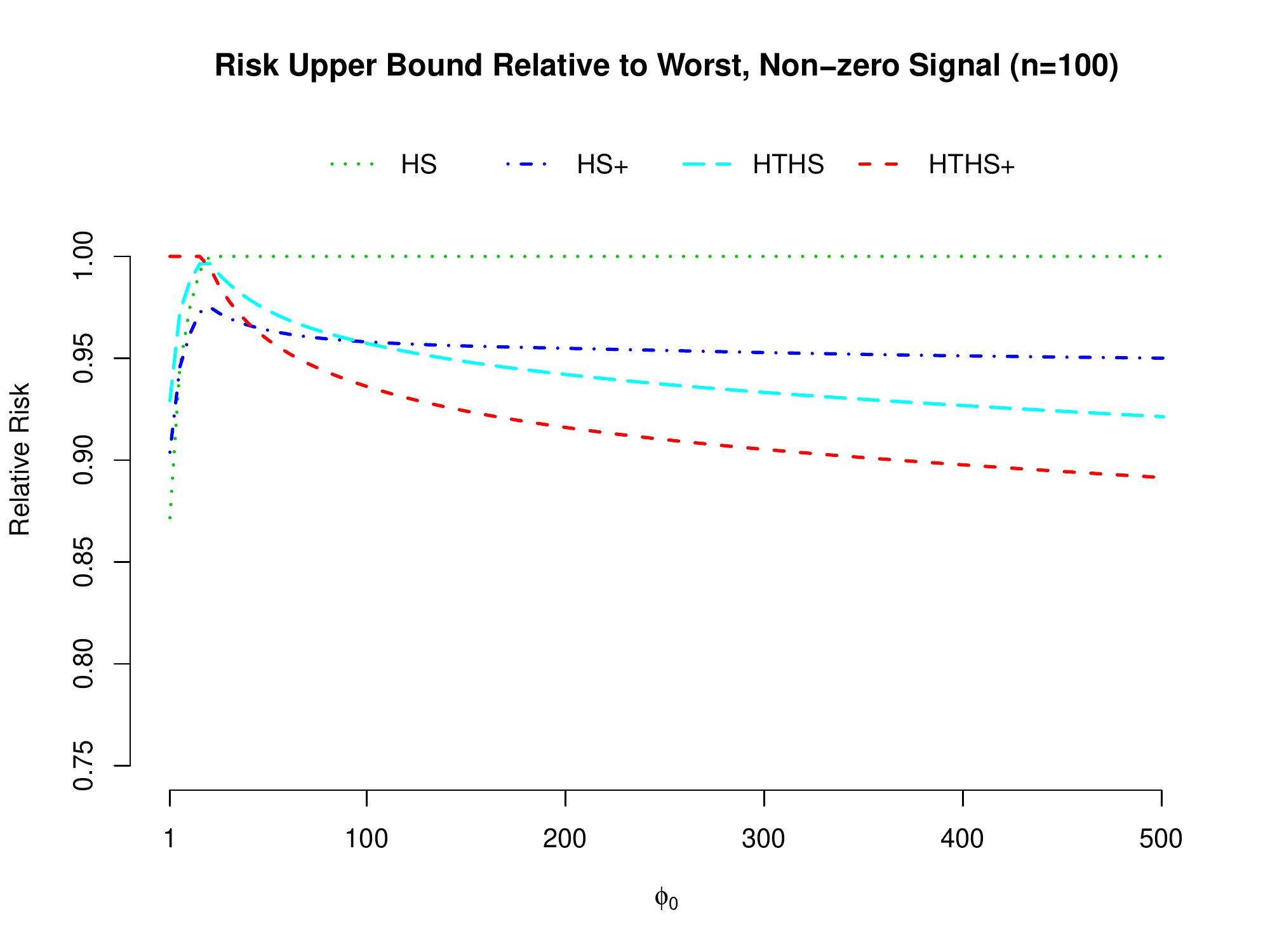} 
\end{minipage}
\caption{\label{fig:KLrisk} KL Risk bounds at the origin and asymptotically for non-zero $\phi_0$.}
\end{figure}

 \begin{thm}
 \label{thm::KL}
 Let $\epsilon=1/n$ and assume that $\phi_0=0$ and that the posterior predictive comes from the HTHS(+) model. Then, the optimal convergence rate $R_n(0)$ is bounded above by
\begin{equation}
\label{eq::riskybusiness}
R_n(0)\le a\frac{\log{n}}{n}+\frac{1}{n}\left(1+\frac{1}{aC_a^0}\right).
\end{equation}
for any $a>0$ where $C_a^0$ is a constant.
\end{thm}
In \cite{bhadra2017horseshoe+}, the authors give the optimal convergence rates of the HS(+) at the true $\phi_0$ as
\begin{equation}
\label{eq::HS_KL}
\begin{array}{c}
R^{HS}_n(0)\le \frac{\log{n}}{2n}-\frac{\log{\log{n}}}{n}+\frac{\text{constant}}{n}\\
R^{HS+}_n(0)\le \frac{\log{n}}{2n}-\frac{2\log{\log{n}}}{n}+\frac{\text{constant}}{n}.
\end{array}
\end{equation}
The second term on the right hand side of \eqref{eq::HS_KL} shows where the improvement is made for HS+ over HS.  Because there is no change to the leading term, one could argue that there is not a significant change in risk made by the HS+ prior. However, this upper bound is somewhat crude and we show the direct computation of the integral in Figure \ref{fig:KLrisk}. The HTHS(+) provides dramatic improvement in the leading term. The shrinking of the first term on the right hand side of \eqref{eq::riskybusiness}, for any $a\in(0,0.5)$, shows that the HTHS(+) attain a   better K-L risk bound and more efficient on suppressing the noise. In Figure \ref{fig:KLrisk}, we also show the relative risks for non-zero $\phi_0$. The differences in efficiency for large signals is not really seen until the signal size is quite large. This shows that the gains in efficiency are really due to better modeling of the true zero signals, especially in ultra-sparse settings. This efficiency feature will be supported by the simulation study in the next section, especially under the ultra-sparse condition. 

\section{Extension}
\begin{figure}[H]
\centering
\begin{minipage}{.5\textwidth}
  \centering
  \includegraphics[width=3in]{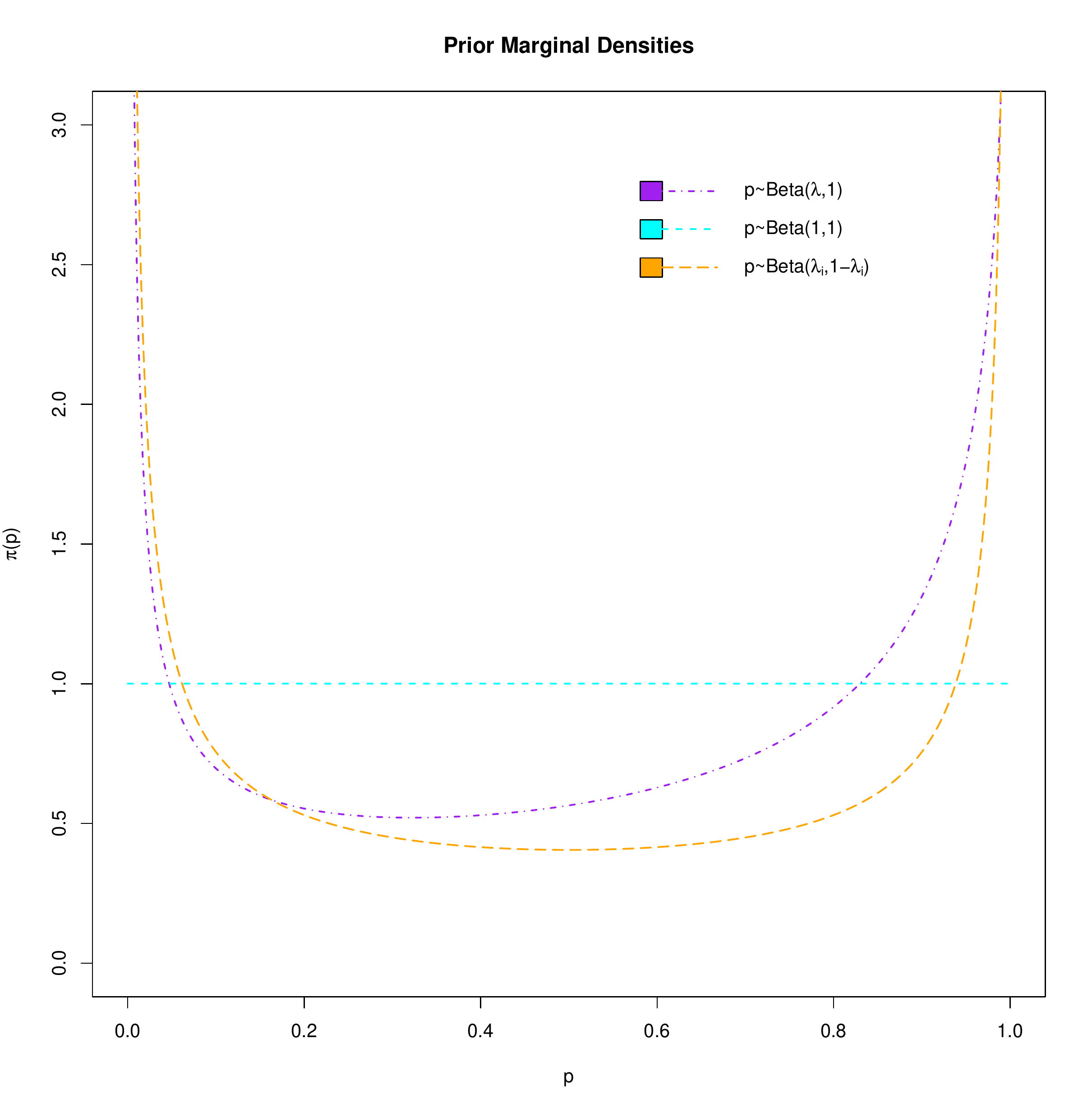}
  \caption{Different prior densities of $p$}
   \label{fig::p_margins}
\end{minipage}%
\begin{minipage}{.5\textwidth}
  \centering
  \includegraphics[width=3in]{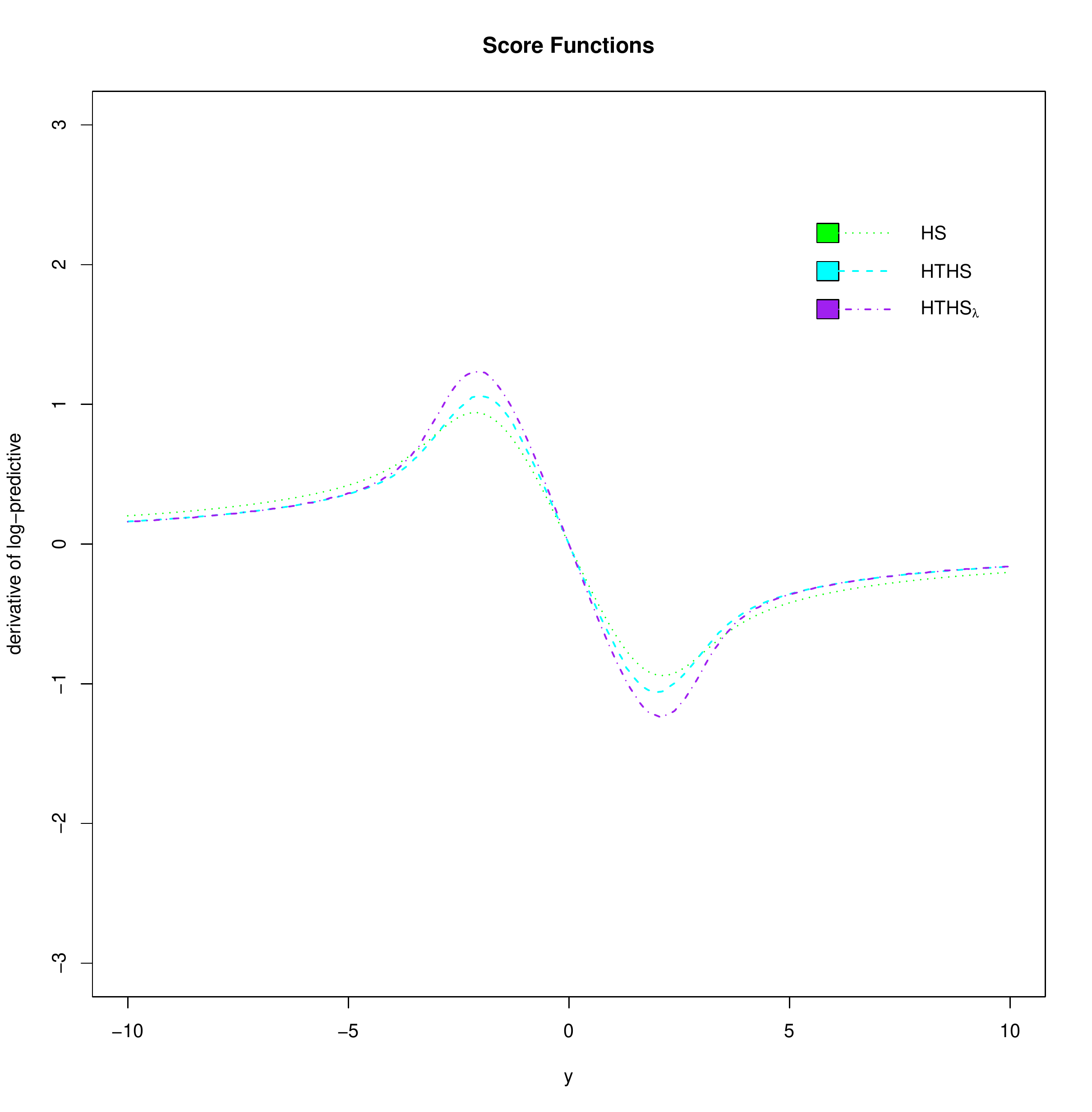}
   \caption{Log-predictive densities}
   \label{fig::score}
\end{minipage}
\end{figure}

Originally, the local decision parameter $p_i$ is assigned with a Beta$(1,1)$ prior distribution, indicating  no particular favors to signals or outliers. However, if the prior knowledge implies a possible ultra-sparse situation, a variant of the HTHS(+) with even stronger suppression on the noises is desired. This motivation is reflected on an asymmetrical $\pi\lrp{p_i}$, further resulting an asymmetrical prior marginal density of the shrinkage profile parameter $\tau_i$. We propose the following prior hierarchy on $p_i$, 
\begin{equation}
\label{prior::asym_p}
\begin{array}{ccc}
p_i|\lambda_i&\sim&\text{Beta}\lrp{\lambda_i,1}\\
\lambda_i|\xi_i&\sim&\text{Gamma}\lrp{1,\xi_i}\\
\xi_i&\sim&\text{Gamma}\lrp{1,1},
\end{array}
\end{equation}
and denote this variant of the HTHS as $\text{HTHS}_\lambda$.
Depending on $\lambda_i$, the prior density of $p_i$ can either have an asymptote at the origin or stack more probability mass near the unit. Unlike $\pi(p_i)$ associated with the HTHS, which is a straight line on the unit interval,  the asymmetry of $\pi(p_i)$ associated with $\text{HTHS}_\lambda$ further reinforces the decision on treating the data points as noises or signals, especially offering extra shrinkages.   A keen reader could notice the similarity between \eqref{prior::asym_p} and the HS or even HTHS hierarchical structure, which is also showing in Fig. \ref{fig::p_margins}.  Indeed, $\pi\lrp{p_i}$ can be easily turned into the HS or the HTHS prior, depending on the situations, but we afraid that including more parameters is more than the model needing now. Notice that $p_i$ is a constant equal to $\frac{1}{2}$ in the HS prior. All the marginal densities of $p_i$ are showed in Figure \ref{fig::p_margins}.\\

 Usually, when the signals are spares, the estimation risk comes from two sources: one from under-shrinking noise, and the other one from over-shrinking the true signals. If the situation is ultra-sparse, this asymmetry structure should provide extra shrinkage, reducing the overall estimation risk. This is exactly showing in Figure \ref{fig::score}, which presents the log-predictive functions of the HS, the HTHS, and the $\text{HTHS}_\lambda$. First, notice that the log-predictive density of the $\text{HTHS}_\lambda$ has  almost identical tails with the HTHS, outperforming the HS on the signal robustness. Further more, introducing $\lambda_i$ makes the $\text{HTHS}_\lambda$ having a smaller risk than the HTHS near the origin, meaning more suitable for coping with the ultra-sparse situation. 
\section{Simulation}
In this section, we show a simple simulation result from the posterior estimation of the HTHS(+), comparing to the result from the HS(+) to demonstrate our theoretical claim in  section \ref{sec::TP}. While exhaustive simulations can be done, this simulation is indicative of results that are obtained for any data analysis using these four models. The data are simulated from the following true model:
\begin{equation}
y_i\sim\lrp{\mu_T+\phi_i,1};~\phi_i\sim \frac{\eta}{2}\text{Uni}(4,6)+\frac{\eta}{2}\text{Uni}(-6,-4)+\lrp{1-\eta}\delta_0\nonumber,
\end{equation}
where $\eta\in\lrc{0.01,0.05,0.2}$ is an overall sparsity parameter and $\delta_0$ is a point mass at zero. All five Bayesian locally adaptive models, the HS(+), the HTHS(+), and the $\text{HTHS}_\lambda$ are included in this simulation to have a straight comparison. They are assigned with same prior distributions for the global parameters to control the discrepancy from any unwanted source. A typical Gibbs with M-H steps is implemented for the posterior distribution samplings with   burning and   thinning procedure.
 We consider two measurements to represent the performance of each model, the mean absolute error and the distance to the oracle M.L.E. estimator. The oracle M.L.E. of $\phi_i$ is defined as $\hat{\phi_i,Ora.}=y_i-\mu_T$ for those $\phi_i$ that are non-zero and $0$ if $\phi_i=0$. The marginal posterior median is adopted for a point estimator under each model due to the multi-modality of each marginal posterior.
\begin{table}[H]

\centering
\begin{threeparttable}
\caption{Bayesian locally adaptive models comparison simulation}
\begin{tabular}{lc|ccc|ccc|ccc|}
\multicolumn{1}{c}{ }&\quad&\multicolumn{3}{|c}{$\eta=$0.2}&\multicolumn{3}{|c}{$\eta=$0.05}&\multicolumn{3}{|c}{$\eta=$0.01}\\
\multicolumn{1}{c}{}&\quad&M.A.E.&\quad&Ora.&M.A.E.&\quad&Ora.&M.A.E.&\quad&Ora. \\
M.L.E.&\quad&321&\quad&262&323&\quad&307&318&\quad&315\\
HS&\quad&198&\quad&146&79&\quad&67&53&\quad&51\\
HS+&\quad&135&\quad&85&65&\quad&53&51&\quad&49\\
HTHS&\quad&134&\quad&84&54&\quad&41&40&\quad&38\\
HTHS+&\quad&95&\quad&46&47&\quad&34&36&\quad&34\\
$\text{HTHS}_\lambda$&\quad&104&\quad&58&33&\quad&22&13&\quad&11
\end{tabular}
\begin{tablenotes}
\item \footnotesize $\eta$, the degree of sparsity; M.A.E., mean absolute error; Ora., distance to the oracle M.L.E. estimator; M.L.E., the maximum likelihood estimator. $n=400$. The numbers are the averages over 20 replicates. 
\end{tablenotes}
\label{tab:simu}
\end{threeparttable}
\end{table}
Table 1 offers a strong evidence for the superiority of adopting the HTHS(+) over the HS(+) when dealing with the sparse situation. Notice that the HS+ indeed has made an improvement over the original HS when the sparsity is moderate, but the two models are really quite close to each other when the situation is ultra-sparsity. In other words, there is no clear advantage for adopting the HS+ over the HS under severe sparsity conditions. As the severity of the sparsity increases, the HTHS(+) makes significant gains over the HS+ in learning noises. It is very interesting to observe that the $\text{HTHS}_\lambda$ has a decisive advantage under ultra-sparse situation comparing to every other models, even though doing worse than the HTHS+ under moderate sparsity. This also shows that the strong asymptote of $\pi\lrp{p}$ in  \eqref{prior::asym_p}  at the origin helps the large signals to escape the extra shrinkage introducing by the asymmetry, preserving the robustness.  
\section{Conclusion}
In this article, we propose an innovative, new Bayesian locally adaptive model that extends the hierarchical structures of the HS and the HS+ priors by adding a local decision parameter. Integrating out the hierarchical parameters yields a log-Cauchy prior distribution for the local adaptive parameter $\gamma_i$ under the HTHS(+) model.  We prove several theorems to theoretical justify the better shrinkage  properties of the HTHS(+). These theoretical properties are exhibited in a simulation study. The next step in this research is to introduce asymmetry into the prior for the local shrinkage profile, which can provide even more dramatic risk gains for ultra-sparse signals.
 \bibliographystyle{apalike}
\bibliography{model_selection.bib} 
\end{document}